\documentclass[a4paper,12pt]{article}

\usepackage{a4}
\usepackage[centertags]{amsmath}
\usepackage{amsfonts}
\usepackage{amssymb}
\usepackage{amsthm}
\usepackage{newlfont}

\begin{document}

\newcommand{\To}{\longrightarrow}

\newtheorem{theorem}{Theorem}[section]
\newtheorem{lemma}[theorem]{Lemma}
\newtheorem{proposition}[theorem]{Proposition}

\theoremstyle{definition}
\newtheorem{remark}[theorem]{Remark}
\newtheorem{remarks}[theorem]{Remarks}

\setlength{\parskip}{0.5cm} \setlength{\unitlength}{0.4cm}
\oddsidemargin 1cm \evensidemargin 1cm \marginparwidth 1in
\marginparsep 7pt \headheight 10pt \headsep .4cm \topmargin 0cm
\textwidth 14.5cm \textheight 23.5cm


\date{}

\author{\textbf{G\"{u}lay Kaya}}
\title{On The Moduli of Surfaces Admitting Genus Two
Fibrations Over Elliptic Curves}
\maketitle

\begin{abstract}

 In this paper, we study the structure, deformations and the moduli
spaces of complex projective surfaces admitting genus two
fibrations over elliptic curves. We observe that, a surface
admitting a smooth fibration as above is elliptic and we employ
results on the moduli of polarized elliptic surfaces, to construct
moduli spaces of these smooth fibrations. In the case of nonsmooth
fibrations, we relate the moduli spaces to the Hurwitz schemes
${\cal H}(1,X(d),n)$ of morphisms of degree $n$ from elliptic
curves to the modular curve $X(d)$, $d\geq 3$. Ultimately, we show
that the moduli spaces in the nonsmooth case are fiber spaces over
the affine line ${\mathbb A}^1$ with fibers determined by the
components of ${\cal H}(1,X(d),n)$.

\end{abstract}


\section{Introduction}

The aim of this paper is to work out the structure, deformations
and the moduli spaces of complex projective surfaces admitting
genus two fibrations over elliptic curves.

In the literature, the cases of albanese fibrations with fiber
genus two over arbitrary base curves and nonalbanese fibrations
over curves of genus $g\geq 2$ have been studied extensively
(\cite{seiler3}, \cite{seiler4} for the former type and
\cite{K-O.pre}, \cite{O-T.2002}, \cite{onsiper.2000} for the
latter). We aim at complementing these results by examining the
case of fibrations with irregularity $q(S)=2$ over elliptic
curves. These fibrations are of nonalbanese type and have Kodaira
dimension $\kappa (S)=1$ (respectively $2$) in case the given
fibration is smooth (respectively non-smooth).

In the smooth case, we fix the type of our surfaces as $\tau =
(1,0;2,2,2;2,1)$, according to the generalized definition of
polarized elliptic surfaces (\cite{seiler2}, p.210). Then we
consider the moduli problem for surfaces $S$ of type $\tau$ as
given. To get a more natural description of this moduli problem we
observe the following. Roughly, forgetting $\tau$, the moduli of
the surfaces we consider is closely related to the moduli of
isogenies of degree two of elliptic curves (the base curves) and
the moduli of smooth genus two curves $C$ admitting an elliptic
subcover $C \to E$ of degree two. The functor corresponding to the
first moduli is coarsely represented by affine modular curve
$Y_0(2)$. The functor corresponding to the latter one is coarsely
represented by an open subscheme $ {\cal H}$ of ${\cal
A}_{2,1}(X(2)\times X(2))/SL_2({\mathbb Z})$ (\cite{kani.1994},
p.210). We obtain

\begin{proposition}\label{prop.1.1}

\begin{itemize}
    \item [(i)] The corresponding functor is coarsely represented
    \ by \ an \ irreducible scheme \ ${\cal M}$ \ of \ dimension $3$ .
    \item [(ii)] There exists a natural surjective morphism
    $\phi: {\cal M} \to {\cal H}$.
\end{itemize}

\end{proposition}

In the nonsmooth case, even though our main interest is in
fibrations over elliptic curves, we will also prove the algebraic
version of the main result in \cite{O-T.2002} describing the
structure of the moduli spaces ${\cal M}(g, K^2, \chi)$ of
surfaces fibered over curves of genus $g \geq 2$. We note that the
result for base curves of genus $g\geq 2$ is stronger than the
result in the case of $g=1$. This is due to the fact that we can
not prove Lemma \ref{lem.2.7} in full strength when $g=1$.

\begin{theorem}\label{thm.1.2}
Let $K^2, ~\chi$ and $g \ge 2$ be given and let ${\cal
H}(g,X(d),n)$ be the Hurwitz scheme of morphisms of degree $n$
from curves of genus $g$ onto $X(d)$. Then we have morphisms $\Phi
: {\cal M}(g, K^2, \chi) \rightarrow {\mathbb A}^1$ and $\Psi_{E'}
: {\cal H}(g,X(d),n) \rightarrow {\cal M}(g, K^2, \chi)$ for any
fixed elliptic curve $E'$ such that \vspace{-.4cm}
\begin{itemize}
    \item [(i)] $\Psi_{E'}$ establishes a one-to-one correspondence
    between the components ${\cal H}_i$ of \\ ${\cal H}(g,X(d),n)$
    and the components ${\cal M}_i$ of ${\cal M}(g, K^2, \chi)$,
    \item [(ii)] $\Phi : {\cal M}_i \rightarrow {\mathbb A}^1$ is
    a fibration with ~$\Psi_{E'}({\cal H}_i)$ ~as the fiber over $[E']
    \in {\mathbb A}^1$.
\end{itemize}
\end{theorem}

\begin{theorem}\label{thm.1.3}
Let ${\cal M}_i$ be a connected component of ${\cal M}
(1,K^2,\chi)$. Then we have a morphism $\Phi : {\cal M}_i \to
{\mathbb A}^1$ $($given on closed points by $[X]\mapsto [E']$ if
$X$ is of type $(E',d))$ such that the fiber over $[E']\in
{\mathbb A}^1$ is a disjoint union
$$\displaystyle\bigsqcup_j~ \Psi_{E'}({\cal H}(1,X(d),n)_j).$$
\end{theorem}

We work over the complex numbers $\mathbb{C}$ and use the
following standard notation:

$S$ is a smooth projective surface.

$c_1(S),~c_2(S)$ denote the first and the second Chern classes of
$S$, respectively.

$\kappa(S),~q(S)$ are the Kodaira dimension and the irregularity
of $S$, respectively.

$K(S),~\chi(S)$ are the canonical class and the holomorphic Euler
characteristic of $S$.

For fixed $K^2$ and $\chi$, ${\cal{M}}(g,K^2,\chi)$ is the moduli
space of surfaces of general type admitting genus two fibration
with irregularity $q=g+1$ and slope $\lambda$ which satisfies the
slope formula $K^2=\lambda \chi + (8-\lambda)(g-1)$.


\section{Structure and Deformations of Genus Two \ Fibrations}

First we discuss the case of smooth fibrations.

\begin{lemma} \label{cor.2.1}
Let $\pi : S \to E$ be a smooth genus two fibration over an
elliptic curve $E$ with $q(S)=2$. Then $S$ admits an elliptic
fibration with two double fibers of the form $2.E$. All other
fibers are smooth and are isomorphic to $E'$ $($the double cover
of $E$ corresponding to the monodromy representation arising from
$\pi)$.
\end{lemma}

\begin{proof}
Let $F$ be the general fiber of the fibration $\pi : S \to E$. Since 
$\pi : S \to E$ is a smooth genus two fibration, $\pi$ is
iso
 with monodromy group $G\subset Aut(F)$ which is cyclic
of order two (\cite{xiao}, Proposition 2.12, p.30) and $F \to
E''=F/G$ ramifies precisely over two points $p_1, p_2 \in E''$ (by
Riemann-Hurwitz formula). Hence, the composite map $E'\times F \to
F \to E''$ induces a natural fibration $S \to E''$ with generic
fiber $E'$ and with two double fibers of the form $2.(E'/G)= 2.E$
over $p_1, p_2$. This proves the Lemma.
\end{proof}

\begin{lemma} \label{lem.2.2}
Let $\psi : {\cal S} \to T$ be a deformation of $S$. Then there
exists an elliptic curve $ {\cal E}'' \to T$ and two sections
$s_1, s_2 :T \to  {\cal E}''$ such that ${\cal S} \to T$ factors
through ${\cal E}''$. Furthermore, ${\cal S} \to {\cal E}''$ is
smooth outside $s_1(T)\cup s_2(T)$ and for each $t\in T$ the
restriction of ${\cal S}\to {\cal E}''$ induces an elliptic
fibration ${\cal S}_t\to {\cal E}''_t$ with precisely two double
fibers $($over $s_1(t), s_2(t))$.
\end{lemma}

\begin{proof}
By standard results in deformation theory, we know that for all
$t\in T$, ${\cal S}_t$ is a minimal surface and $\kappa ({\cal
S}_t)=1$. Furthermore, each ${\cal S}_t$ admits an elliptic
fibration exactly of the same type as $S$ (\cite{FM}, Proposition
7.1, p.111) and it follows from (\cite{FM}, Proposition 7.11(iii),
p.128) that there exists an elliptic curve $ {\cal E}'' \to T$
through which ${\cal S} \to T$ factors. Since each surface ${\cal
S}_t$ is an ``elliptic surface of general type" in the terminology
of ~\cite{itaka}~ (i.e., $P(n)\geq 2$ for some $n$), Proposition
10 in \cite{itaka} applies to prove the existence of two sections
$s_1, s_2:T \to {\cal E}''$ (\cite{seiler2}, Lemma 1.9) with the
properties stated in the Lemma.
\end{proof}

\begin{remarks}
\begin{itemize}
    \item [\textbf{1)}] The existence of ~${\cal E}''$ over $T$ can
be proved simply by observing that the fibrations in the family
${\cal S}$ are induced from $m$-th canonical map for $m$
sufficiently large (\cite{seiler2}, p.194). More precisely, we
take ${\mathbb P}(\psi_* \omega^{\otimes m}_{{\cal S}/T})$ over
$T$ and the morphism ${\cal S} \to {\mathbb P}(\psi_*
\omega^{\otimes m}_{{\cal S}/T})$ induced by the homomorphism
$\psi^*(\psi_* \omega^{\otimes m}_{{\cal S}/T}) \to
\omega^{\otimes m}_{{\cal S}/T} \to 0$. ~${\cal E}''$ is the image
of this $m$-th canonical map.
    \item [\textbf{2)}] Lemma \ref{lem.2.2} indicates a relation
between the moduli of smooth genus two fibrations over elliptic
curves with irregularity $q=2$ and the moduli of elliptic
surfaces. In order to apply the results of (\cite{seiler1},
\cite{seiler2}), we need the following observation:
\end{itemize}
\end{remarks}

\begin{lemma} \label{lem.2.3}
An elliptic surface with exactly two double fibers as the singular
fibers, admits a smooth genus two fibration.
\end{lemma}

\begin{proof}
Let $\pi :S \to E''$ be an elliptic fibration over an elliptic
curve $E''$ with two double fibers over $p_1, p_2 \in E''$ and
generic fiber $E'$. Let $\pi _J: B=J(S) \to E''$ be the Jacobian
fibration and $B'=E'\times E'' \to E''$ be the trivial elliptic
fibration. Both, $B$ and $B'$, are elliptic fibrations with
sections. Since the associated $j$-invariants are equal and
constant, we can find isomorphic compatible lifts $\rho$ and $\rho
'$ defined by $B$ and $B'$, respectively (\cite{FM}, p.41). Hence,
the Jacobian surface $J(S)$ of $S$ is trivial, being isomorphic to
$B'=E'\times E''$ (\cite{FM}, Theorem 3.14(ii), p.45). Moreover,
$R^1\pi_* {\cal O}_S \cong R^1(\pi_J)_* {\cal O}_{J(S)}$ is
trivial and so $L=(R^1\pi_* {\cal O}_S)^\vee$ is trivial.
Therefore, $\chi({\cal O}_S)=\deg L =0$ (\cite{FM}, Proposition
3.18, p.48), $p_g(S)=g(E'')=1$ (\cite{FM}, Proposition 3.22(i),
p.49) and so $q(S)=1+p_g(S)-\chi({\cal O}_S)=2$. By the universal
property of albanese varieties there is a morphism $Alb(S)\to
E''$. Hence, $Alb(S)$ is a reducible abelian variety and by the
complete reducibility property of abelian varieties
(\cite{mumford.70}, Theorem 1, p.173) it has a projection to an
elliptic curve $E$, $\pi_1 : Alb(S) \to E$, which restricts to a
nonconstant necessarily \'{e}tale morphism from $E'$. The pullback
$S\times_E E'$ is the product $F\times E'$ where $F$ is
the general fiber of the composite map $S\to Alb(S) \to E$
(\cite{kollar}, E.8.6, p151). $g(F)=2$ since $F$ is a double cover
of $E''$ which is ramified at two points.
\end{proof}

\remark Combining the results in Lemma \ref{lem.2.2}~ and Lemma
\ref{lem.2.3}, we see that a surface $S$ admitting a smooth genus
two fibration deforms only to surfaces of the same type.

Next, we consider the case of nonsmooth fibrations. Let $\pi :
S\to C$ be a nonsmooth genus two fibration with $q(S)=g(C)+1$.
Then there is a unique rational number $\lambda = \lambda(\pi)$,
which is called the slope of $\pi$, such that
$K^2=\lambda\chi+(8-\lambda)(g(C)-1)$, and one has
$2\leq\lambda\leq 7$ for nontrivial nonsmooth fibrations
(\cite{xiao}, p.22). Let $F$ be a smooth fiber of $\pi$ and $J(F)$
its Jacobian. We have a projection $p_F: J(F) \to E$ onto the
fixed part of the associated relative Jacobian. Let $d$ be the
degree of the composite map $F \to J(F) \to E$. Then $d$ is called
the degree associated to $\pi$ and $\pi : S\to C$ is said to be of
type $(E,d)$. We have $\lambda=7-\displaystyle\frac{6}{d}$ \
(\cite{xiao}, Corollaire 2, p.50). Since we will consider only
semistable fibrations, we have $d\geq 3$ (\cite{xiao}, Corollaire,
p.47). As an immediate result of (\cite{xiao}, Th\'{e}or\`{e}me
3.10, p.44) we obtain

\begin{theorem}\label{thm.2.1}
Let $E$ be an elliptic curve, d an integer $\geq 3$. There exists
a genus two fibration of type $(E,d)$ $$ \Phi : S(E,d) \to X(d)$$
over the modular curve $X(d)$ which is universal in the following
sense: any genus two fibration $\pi : S \to C$ with slope $\lambda
=7- \displaystyle \frac{6}{d}$ and with $E$ as the fixed part of
the Jacobian fibration corresponding to $\pi$ $($i.e., $\pi$ is of
type $(E,d))$ is the minimal desingularization of the pullback
$f^*(S(E,d))$ via a surjective holomorphic map $f:C\to X(d)$.
\end{theorem}

\remark Since $g(X(d))\geq3$ for $d\geq 7$, it follows that
fibrations over elliptic curves have $d \in \{3,4,5,6\}$. We
recall that $X(d)\cong {\mathbb P}^1$ for $d=3,4,5$ and $X(6)$ is
the elliptic curve with $j(X(6))=0$.

Given $f:C\to X(d)$, the surface $f^*(S(E,d))$ has singularities
only if $f$ ramifies over some points in the singular locus of $
\Phi : S(E,d) \to X(d)$. A singular fiber of $ \Phi$ is either an
elliptic curve with a single node or two smooth elliptic curves
intersecting transversally at a single point (\cite{xiao}, Lemme
3.11, Th\'{e}or\`{e}me 3.16). Hence, singularities of
$f^*(S(E,d))$ are all type $A_k$ for some $k$ depending on the
singular point.

As a consequence of this observation we see that we can apply
simultaneous desingularization to a family of surfaces obtained
via a family of surjective morphisms onto $X(d)$.

\begin{lemma} \label{lem.2.4}
For a fibration $\pi : S \to C$ over a curve $C$ of genus $\geq 1$
arising from a map $f:C\to X(d)$ of degree $n$ we have $c_2(S)>0$
and $K^2 = {c_1}^2(S)>0$. In particular, since $S$ is minimal, it
is a surface of general type.
\end{lemma}

\begin{proof}
Let $\phi : S(E',d)\to X(d)$ be the corresponding fibration. Then
$c_2(S)=-n\deg(R^1\phi_*{\cal O}_{S(E',d)})>0$. Using the relations
$c_1^2(S)=\lambda\chi(S)+(8-\lambda)(g(C)-1)$ and
$12\chi(S)=c_1^2(S)+c_2(S)$, we have $c_1^2(S)>0$.
\end{proof}

\begin{lemma} \label{lem.2.5}
Let $S_i \rightarrow C_i$, $i = 1,2$ be two fibrations of the same
type $(E, d)$, corresponding to morphisms $f_i : C_i \rightarrow
X(d)$. Then \vspace{-.4cm}
\begin{itemize}
    \item [(i)] $S_i$ have the same invariants $K^2,~\chi$ if and only if
    $\deg(f_1) = \deg(f_2)$,
    \item [(ii)] $S_1$ and $S_2$ are isomorphic as {\it surfaces} if and only
if ~$C_1 = C_2$ and there exist automorphisms $\alpha \in
Aut(C_1), ~\beta \in Aut(X(d))$ such that $f_1 \circ \alpha =
\beta \circ f_2$.
\end{itemize}
\end{lemma}

\begin{proof}
\noindent\textbf{(i)} This is Lemma 1 in \cite{O-T.2002}.

\noindent\textbf{(ii)} That $C_1 = C_2$ follows from the
uniqueness of such a fibration on a given surface
(\cite{O-T.2002}, Lemma 2(i)). Then the rest of the statement is a
consequence of the minimality of the surfaces $S_1$ and $S_2$,
since the given condition is necessary and sufficient for the
surfaces $f_i^*(S(E,d))$ to be birationally equivalent.
\end{proof}

\begin{lemma} \label{lem.2.6}
$S$ admitting a fibration as described over an elliptic curve,
exists if and only if $K^2$ and $\chi$ have the following values:

\begin{tabular}{c c c}
$\lambda$ & \qquad $K^2$ & \qquad $\chi$ \\
$5$ & \qquad $5n$ & \qquad $n$\\
$11/2$ & \qquad $11n$ & \qquad $2n$\\
$29/5$ & \qquad $29n$ & \qquad $5n$\\
$6$ & \qquad $36n$ & \qquad $6n$
\end{tabular}

\vspace{.5cm} \noindent where $n\geq 2$ in the first three rows
and $n\geq 1$ in the last row.
\end{lemma}

\begin{proof}
For $n\geq 2$ and for any elliptic curve $E$ we have morphisms
$E\to {\mathbb P}^1$ of degree $n$ and for any such a map, using
the formulae in the proof of Lemma \ref{lem.2.4} and observing
that ~$-\deg(R^1\phi_*{\cal O}_{S(E',d)})=7,13,31$ for $d=3,4,5$,
respectively, (\cite{xiao}, p.52), we find the values of $K^2$ and
$\chi$ given in the first three rows of the table.

Since $X(6)$ is an elliptic curve, by the same computation, this
time using the existence of isogenies of any order and the fact
that $-\deg(R^1\phi_*{\cal O}_{S(E',6)})=6$ we obtain the last row.
\end{proof}

We will need the following Lemma (\cite{O-T.2002}, Lemma 2).

\begin{lemma} \label{lem.2.7}
Let  $\psi : {\cal S} \rightarrow T$ be a deformation over a
connected base, of a surface $S$ admitting a genus 2 fibration
with slope $\lambda$ over a curve $C$ of genus $g \geq 2$. Then

\begin{itemize}
    \item [(i)] each fiber ${\cal S}_{t}$ of $\psi$ admits such a
fibration ${\cal S}_{t} \rightarrow C_{t}$ which is unique,
    \item [(ii)] the slope ~$\lambda$ is constant on $T$,
    \item [(iii)] the degree of the map $C_{t} \rightarrow X(d)$ inducing
the fiber space ${\cal S}_{t} \rightarrow C_{t}$ is constant.
\end{itemize}
\end{lemma}

In case of elliptic base curves (ii) and (iii) of Lemma
\ref{lem.2.7} remain unchanged when we consider a family ${\cal S}
\to T$ of surfaces having a fibration of the given form. Moreover,
the fibration over any such curve is also unique. However, we do
not know if (i) holds, too.


\section{Moduli Problem of Genus Two Fibrations}

Let $\pi : S\to E$ be a smooth genus two fibration with fiber $F$
and $\pi_1 : S\to E''$ be the corresponding elliptic fibration
with two double fibers. In fact, these double fibers are sections
of $\pi$, say $s_1, s_2$. Consider the divisor $s_1(E)+F$ on $S$.
We have $$(s_1(E)+F)^2=2s_1(E).F=2>0$$ and for any irreducible
curve  $C$ in $S$ $$(s_1(E)+F).C=s_1(E).C+F.C>0.$$ Hence,
$s_1(E)+F$ is an ample divisor on $S$, by Nakai's ampleness
criterion. Let $\eta$ be the numerical equivalence class of the
line bundle corresponding to this ample divisor in $Num(S)$ (the
group of numerical equivalence classes of line bundles on $S$).
Then $d:=\eta^2=(s_1(E)+F)^2=2$ and
$e:=\eta.f=(s_1(E)+F).E'=F.E'=1$ where $f$ is the class of the
general fiber $E'$ of $\pi_1$. Hence, our surfaces are of type
$\tau=(1,0;2,2,2;2,1)$ according to the generalized definition of
polarized elliptic surfaces given in (\cite{seiler2}, p.210).
Moreover, any surface of type $\tau$ is one of our surfaces.

We consider the functor $G_{\tau} : Sch \to Sets$ defined by
$G_{\tau}(T)=$ set of all isomorphism classes of families of
polarized elliptic surfaces of type $\tau$ over $T$. We have

\begin{proposition} \label{prop.3.1}
$G_{\tau}$ is coarsely represented by an irreducible scheme ${\cal
M}$ of dimension 3.
\end{proposition}

\begin{proof} Existence of ${\cal M}$ follows from (\cite{seiler2},
Theorem 2.15, p.211). The proof of this theorem shows that there
is a finite map ${\cal M} \to Y''$, where $Y''$ is an open
subscheme of $Y'=E_{1,0}\times_{{\mathbb A}^1} M_{1,2}$. Here
$E_{1,0}$ denotes the moduli scheme for Weierstrass surfaces with
base genus $g=1$ and $\chi=0$ which exists by \cite{seiler1}, and
$M_{1,2}$ is the  moduli scheme for elliptic curves with two
distinguished points. $E_{1,0}$ splits into a disjoint union of
irreducible subschemes $E_{1,0}^n$ for $n=1,2,3,4,6$
(\cite{seiler1}, p.182) where each $E_{1,0}^n$ represents the
subfunctor corresponding to Weierstrass surfaces for which the
order of the module $L=(R^1p_*{\cal O})^\vee$ is $n$.

Let $S$ be an elliptic surface of type $\tau$. Since all fibers of
the elliptic fibration on $S$  are irreducible, the Weierstrass
fibration associated to $S$ is the Jacobian fibration of $S$
(\cite{seiler2}, p.191). In the proof of Lemma \ref{lem.2.3}, we
have seen that the Jacobian of such a surface is a trivial product
of two elliptic curves. So the relevant part of $E_{1,0}$ is
$E_{1,0}^1$ which corresponds to trivial $L$. Hence, $Y''$ is an
open subscheme of $E_{1,0}^1\times_{{\mathbb A}^1} M_{1,2}$. By
(\cite{seiler1}, Lemma 10, p.182) $E_{1,0}^1\cong {\mathbb A}^2$.
Therefore, $\dim({\cal M})=\dim(E_{1,0}^1\times_{{\mathbb A}^1}
M_{1,2})=3$.
\end{proof}

In the preceding section we have observed that the moduli of the
surfaces we consider is closely related to the moduli of isogenies
of degree two of elliptic curves (the base curves) and the moduli
of smooth genus two curves $C$ admitting an elliptic subcover $C
\to E$ of degree two. The functor ${\cal Y}_0$ corresponding to
the first moduli is coarsely represented by affine modular curve
$Y_0(2)$. As for the latter, we have the affine surface ${\cal A}
_{2,1} = (X(2)\times X(2))/SL_2({\mathbb Z})$ (\cite{kani.1994},
p.210) which coarsely represents the functor associated to the
triplets $\left\{ (A , \Theta, E) \right\}$ where $A$ is an
abelian surface, $\Theta$ is a principal polarization, $E$ is an
elliptic subgroup of $A$ and $\deg(\Theta|_E)=2$. Let $C$ be a
curve of genus two with Jacobian $J_C$ and  canonical polarization
$\Theta$. Then there is a bijective correspondence between the set
of isomorphism classes of (minimal) elliptic subcovers $f:C\to E$
of degree $\deg(f)=2$ and the set of elliptic subgroups $E\leq J_C$
of $J_C$ of degree $\deg_{\Theta}(E)=2$ (\cite{kani.1994}, Theorem
1.9, p.202). Therefore, the functor ${\cal M}'_2$ of isomorphism
classes of pairs $({\cal C}, {\cal E})$ of (relative) smooth
curves of genus two and elliptic subcovers $({\cal C}\to {\cal
E})$ of degree two is coarsely represented by the open subscheme $
{\cal H}= \Phi ^{-1} (t({\cal M}_2))$ of ${\cal A}_{2,1}$, where
$t:{\cal M}_2 \to {\cal A}_2$ is the Torelli map which associates
to a curve its canonically polarized Jacobian and $\Phi : {\cal
A}_{2,1} \to {\cal A}_2$ is the map which forgets $E$ in the
triplets.

\begin{proposition} \label{prop.3.2}
There exists a natural surjective morphism $\phi: {\cal M} \to
{\cal H}$.
\end{proposition}

\begin{proof} Consider an object in $G_{\tau}(T)$ for some $T$;
i.e. a family ${\cal S}\to T$ which factors over ${\cal E}''$
(Lemma \ref{lem.2.2}). Let $f_1,f_2\ \in {\cal O}_{{\cal E}''/T}$
such that $(f_i)=s_i(T),\ i=1,2$. Then the natural injection
$f:{\cal O}_{{\cal E}''/T} \to {\cal O}_{{\cal E}''/T} [\sqrt{f_1
f_2}~]$ gives a double cover $(f):{\cal F}/T \to {\cal E}''/T$
over $T$ ramified along $(f_1)\cup(f_2) = s_1(T)\cup s_2(T)$,
where ${\cal F}=Spec~({\cal O}_{{\cal E}''}[t]/(t^2-f_1 f_2))$.
Moreover, since for any $t\in T$, ~${\cal F}_t \to {\cal E}''_t$
is a double cover ramified at two points, we have $g({\cal
F}_t)=2$ by Riemann-Hurwitz formula. Hence, ${\cal F}_t$ is a
smooth genus two curve. Hence, the pair ${({\cal F}, {\cal E}'')}$
corresponds to a point in ${\cal H}(T)$. By functoriality of this
construction, we obtain natural morphism $\phi: {\cal M} \to {\cal
H}$.
\end{proof}

Next we consider moduli of surfaces with nonsmooth genus two
fibrations of nonalbanese type. In Section 2 we have seen that a
surface of type $(E',d)$ is the desingularization of
$f^*(S(E',d))$ for some morphism $f:C \to X(d)$. Hence, such a
surface $S$ can be deformed in two ways; we can deform $E'$ to
other elliptic curves and we can deform the map $f$. Therefore, in
describing the moduli spaces of such surfaces under consideration,
we need to clarify the relation of these spaces to the Hurwitz
spaces ${\cal H}(g,X(d),n)$ of morphisms of degree $n$ from curves
of genus $g$ to the modular curve $X(d)$.

\begin{theorem} \label{thm.3.1}
Let $K^2, ~\chi$ and $g \ge 2$ be given and let ${\cal
H}(g,X(d),n)$ be the Hurwitz scheme of morphisms of degree $n$
from curves of genus $g$ onto $X(d)$. Then we have morphisms $\Phi
: {\cal M}(g, K^2, \chi) \rightarrow {\mathbb A}^1$ and $\Psi_{E'}
: {\cal H}(g,X(d),n) \rightarrow {\cal M}(g, K^2, \chi)$ for any
elliptic curve $E'$ such that
\begin{itemize}
    \item [(i)] $\Psi_{E'}$ establishes a one-to-one correspondence
    between the components ${\cal H}_i$ of ${\cal H}(g,X(d),n)$
    and the components ${\cal M}_i$ of ${\cal M}(g, K^2, \chi)$,
    \item [(ii)] $\Phi : {\cal M}_i \rightarrow {\mathbb A}^1$ is
    a fibration with ~$\Psi_{E'}({\cal H}_i)$ ~as the fiber over $[E']
    \in {\mathbb A}^1$.
\end{itemize}
\end{theorem}

This result is a consequence of Lemma \ref{lem.3.1} and Lemma
\ref{lem.3.2}.

\begin{lemma} \label{lem.3.1}
There exists a morphism $\Phi :{\cal M}(g, K^2, \chi) \rightarrow
{\mathbb A}^1$ which maps the class $[S] \in {\cal M}(g, K^2,
\chi)$ to the class $[E] \in {\mathbb A}^1$ of the elliptic curve
associated to the fibration on $S$. $\Phi$ is surjective on each
component of ${\cal M}(g, K^2, \chi)$.

\end{lemma}

\begin{proof}
Let $M(g,~\lambda): Sch/{{\mathbb C}} \rightarrow Sets$ be the
functor defined by ~$M(g,~\lambda)(T)$ = isomorphism classes of
families of surfaces over T admitting genus 2 fibrations over
curves of genus $g$, with slope $\lambda$. To prove the lemma, it
suffices to construct a morphism of functors ~$M(g,~\lambda)
\rightarrow h_{{\mathbb A}^1}$ as described in the lemma. This, on
the other hand, follows once we prove that for any ~$T \in Sch
/{\mathbb C}$~ and for a given family of surfaces ~${\cal S}
\rightarrow T$, the map ~$T \rightarrow {\mathbb A}^1$ ~defined by
~$t \mapsto [E_t]$ ~where $[E_t]$ is the fixed part of the
jacobian fibration on ~${\cal S}_t$, is a morphism.

This last claim being local over the base, we assume that ~${\cal
S} \rightarrow T$ is projective and we consider the relative
albanese morphism ~$\alpha : {\cal S} \rightarrow Alb_{{\cal
S}/T}$; the image ~${\cal E} = \alpha({\cal S})$~ is a family of
smooth isotrivial elliptic surfaces over $T$ and the base of the
fibration on ~${\cal E}_t$ is $C_t$ = the base of the fibration on
${\cal S}_t$ . It is well known that for such a family of elliptic
surfaces, the base curves glue to give a relative curve ${\cal C}$
and the map ~${\cal E} \rightarrow T$ factors over $\cal C$. Since
the fibres of ${\cal E}_t \rightarrow C_t$ are constant, the
morphism ~${\cal C} \rightarrow {\mathbb A}^1$ corresponding to
the elliptic curves ${\cal E}/{\cal C}$ coincides with the map ~$T
\rightarrow {\mathbb A}^1$ defined above, which completes the
proof of the claim.

To~~prove~~the~~surjectivity~~of ~$\Phi$,~
we~~take~~any~~connected~~component~~of ${\cal M}(g,K^2,\chi)$ and
a surface $S$ of type $(E, d)$ corresponding to a point in this
component. We let $f : C \rightarrow X(d)$ be the map inducing the
fibration on $S$. For any family of elliptic curves ~${\cal E}
\rightarrow T$, we have a genus two curve ~${\cal F} \rightarrow
H_{{\cal E}/T, d,-1}$, where ~$H_{{\cal E}/T, d,-1}$ is an open
subscheme of ~$X(d) \times_{{\mathbb C}}T$, universal for
normalized genus two covers ([1], Definition on p.13) of degree
$d$ of ~${\cal E}/T$ ([1], Thm. 1.1). From $(f,id) : C \times T
\rightarrow X(d) \times_{{\mathbb C}}T$ we obtain a $T$-morphism
~$F : U \rightarrow H_{{\cal E}/T, d,-1}$ where $U$ is an open
subscheme of $C \times T$. Completing the family of genus two
curves ~$F^*({\cal F})/U$ ~to a family over ~$C \times T$, and
then applying simultaneous desingularization we get a family of
smooth surfaces ~${\cal S} \rightarrow T'$ where $T' \rightarrow
T$~ is a finite Galois base extension. Since ~${\cal S}$ contains
$S$ as one of the fibers, its moduli lies in the same component of
~${\cal M}(g, K^2, \chi)$~ as the modulus of $S$. For an arbitrary
elliptic curve ~$E'$, choosing ~${\cal E} \rightarrow T$ as a
deformation of $E$ to $E'$, we see that ~$\Phi$~ restricted to
this moduli has $[E'] \in {\mathbb A}^1$ in its image. This
completes the proof of the lemma.
\end{proof}

Let ${\cal C} \rightarrow T$ be a family of smooth curves of genus
$g$ and let $F : {\cal C} \rightarrow X(d)\times T$ be a family of
morphisms of degree $n$. For a fixed elliptic curve $E'$, applying
simultaneous desingularization to the family of surfaces
$F^*(S(E', d))$  we obtain a family of fibered surfaces ${\cal S}
\rightarrow T'$ over a Galois extension $T' \rightarrow T$ with
group $G$, which defines a morphism $\alpha : T' \rightarrow {\cal
M}(g, K^2, \chi)$. $\alpha$, being $G$-invariant, descends to a
morphism $T \rightarrow {\cal M}(g, K^2, \chi)$. Clearly, this
construction is functorial and by the defining property of coarse
moduli spaces we get a morphism ~$\Psi_{E'} : {\cal H}(g,X(d),n)
\rightarrow {\cal M}(g, K^2, \chi)$. To a given connected
component ${\cal H}_i$ of ${\cal H}(g,X(d),n)$ we assign the
component ${\cal M}_i$ of ${\cal M}(g, K^2, \chi)$ which contains
~$\Psi_{E'}({\cal H}_i)$.

\begin{lemma} \label{lem.3.2}
The above assignment induces a one-to-one correspondence between
the connected components of ${\cal M}(g, K^2, \chi)$ and those of
${\cal H}(g,X(d),n)$. Morever, we have $\Psi_{E'}({\cal H}_i)=
\Phi ^{-1} \vert _{{\cal M}_i}([E'])$.
\end{lemma}

\begin{proof}
Since by Lemma \ref{lem.3.1}, each component ${\cal M}_i$ of
${\cal M}(g, K^2, \chi)$ contains the modulus of a surface of type
$(E', \lambda)$, it suffices to check that in each ${\cal M}_i$ we
have the image under $\Psi_{E'}$ of a unique component of ${\cal
H}(g,X(d),n)$.

Let ${\cal M}_i$ be  the component of ~${\cal M}(g, K^2, \chi)$
~which contains ~$\Psi_{E'}({\cal H}_i)$. Fix $[S_1] \in
\Psi_{E'}({\cal H}_i)$ ~and let ~$[S_2] \in {\cal M}_i$ be an
arbitrary point and let ~$C_i, ~j=1,2$~ be the base curves of the
corresponding fibrations. Then, the surfaces $S_1$ and $S_2$
deform to each other. Since deformations of the surfaces under
consideration are induced from deformations of the fibrations
(proof of Lemma \ref{lem.3.1}), it follows that ~$f_1 : C_1
\rightarrow X(d)$ deforms to a morphism ~${\overline f}_2 : C_2
\rightarrow X(d)$. Therefore, ~$f_1,~{\overline f}_2$ belong to
~${\cal H}(g,X(d),n)_i$. On the other hand, by (Lemma
\ref{lem.2.5} (ii)), $f_2$ ~and ~${\overline f}_2$ satisfy a
relation of the form ~$f_2 \circ \alpha = \beta \circ {\overline
f}_2$ for some ~$\alpha \in Aut(C_1), ~\beta \in Aut(X(d))$.
Therefore, $f_2$~and ${\overline f}_2$, hence, $f_1$ and $f_2$ lie
in ~${\cal H}_i$. This proves the first part of the lemma. The
second statement is obvious.
\end{proof}

In the case of elliptic base curves we can not prove that a given
deformation of our surfaces ${\cal S} \to T$ arises from the
deformation of the associated maps $f_t: C_t \to X(d),\ t\in T$.
Therefore, by exactly the same proof we obtain the following
weaker result:

\begin{theorem}\label{thm.3.2}
Let ${\cal M}_i$ be a connected component of ${\cal M}
(1,K^2,\chi)$. Then we have a morphism $\Phi : {\cal M}_i \to
{\mathbb A}^1$ $($given on closed points by $[X]\to [E']$ if $X$
is of type $(E',d))$ such that the fiber over $[E']\in {\mathbb
A}^1$ is a disjoint union
$$\displaystyle\bigsqcup_j~ \Psi_{E'}({\cal H}(1,X(d),n)_j).$$
\end{theorem}

\begin{remarks}
\begin{itemize}
    \item [\textbf{1)}] When $\lambda=6$, one can prove that ${\cal
M}(g,K^2,\chi)=\displaystyle\bigsqcup_{i=1}^N~{\mathbb A}^1_i$
where $N$ is the number of distinct \'etale covers of degree $n$
of the elliptic curve $X(6)$ (\cite{karadogan}, Theorem 2.3).
    \item [\textbf{2)}] Another shortcoming of the result in case of base genus
$g=1$ is that we do not know if each ${\cal M}_i$ is a connected
component of the corresponding moduli space ${\cal M}_{K^2,\chi}$
of surfaces of general type.
\end{itemize}
\end{remarks}

\vskip .5cm \noindent G\"{u}lay Kaya\\
Department of Mathematics\\
Galatasaray University\\
Ciragan Cad. No:36 34357
Ortakoy, Istanbul, Turkey

\label{`@lastpage'}

\end{document}